\documentclass[12pt]{article}
\usepackage{amsmath}
\def\C{\mathbf{C}}
\def\Z{\mathbf{Z}}
\def\mod{\mathrm{mod}\, }
\def\R{\mathbf{R}}
\title{Exceptional solutions to the Painlev\'e VI equation}
\author{Alexandre Eremenko\thanks{Supported by NSF grant DMS--1361836},
Andrei Gabrielov\thanks{Supported by NSF grant DMS-1161629}\; and Aimo Hinkkanen}
\begin{document}

\maketitle
\begin{abstract}
We find all solutions of the Painlev\'e VI equations with the property
that they have no zeros, no poles, no $1$-points and no fixed points.

2010 AMS Class.: 34M55, 33E05. Keywords: Painlev\'e equations,
elliptic functions,
modular group.
\end{abstract}

Painlev\'e VI is the following second order ODE:
\begin{equation}\label{1}
\begin{aligned}
\frac{d^2y}{dt^2}&=\frac{1}{2}\left(\frac{1}{y}+\frac{1}{y-1}+\frac{1}{y-t}
\right)\left(\frac{dy}{dt}\right)^2-
\left(\frac{1}{t}+\frac{1}{t-1}+\frac{1}{y-t}\right)\frac{dy}{dt}\\
&+\frac{y(y-1)(y-t)}{t^2(t-1)^2}\left[\alpha+\beta\frac{t}{y^2}+\gamma
\frac{t-1}{(y-1)^2}+\delta\frac{t(t-1)}{(y-t)^2}\right],
\end{aligned}
\end{equation}
where $(\alpha,\beta,\gamma,\delta)$ are complex parameters.

It is known that each solution has a meromorphic continuation
along every curve in $D=\C\backslash\{ 0,1\}$, see, for example, \cite{H}.
A solution $y(t)$ is called {\em exceptional} if
$y(t)\not\in \{ 0,1,\infty,t\}$ for all $t\in D$ (and for all branches of $y$).
In \cite{Lin} such solutions are called ``smooth''.

An interesting problem is to classify all exceptional solutions.

When
$$(\alpha,\beta,\gamma,\delta)=(0,0,0,1/2),$$
equation (\ref{1}) was studied by Picard \cite{Picard} 16 years before
Painlev\'e, Gambier and R. Fuchs discovered it. Picard found all solutions
for this case
(see below). All Picard's solutions are exceptional. 
The following two results are known.
\vspace{.1in}

When
$$(\alpha,\beta,\gamma,\delta)=(1/8,-1/8,1/8,3/8),$$ there are exactly three
exceptional solutions \cite{Lin}. 
\vspace{.1in}

Local solutions are considered the same if they
are obtained by an analytic continuation from each other.
\vspace{.1in}

When
$$(\alpha,\beta,\gamma,\delta)=(9/8,-1/8,1/8,3/8)$$
there is exactly one exceptional solution defined by the equation
\begin{equation}\label{lin}
3y^4-4ty^3-4y^3+6ty^2-t^2=0.
\end{equation}
This was recently found in \cite{Lin2}.
\vspace{.1in}

We give a simple proof of these results. Moreover,
we determine all values of parameters for which exceptional solutions
exist, find their number for such values of parameters,
and write down explicit representations of these solutions.

It will be convenient to work with the {\em elliptic form} of Painlev\'e VI discovered by
Picard.
Consider the lattice $\Lambda_\tau=\{ m+n\tau:m,n\in\Z\}$, where $\tau$ is
in the upper half-plane $H$.
The Weierstrass function $\wp(z|\tau)$ is the solution of the
differential equation
$$(\wp')^2=4(\wp-e_1)(\wp-e_2)(\wp-e_3),$$
with the initial condition $\wp(0)=\infty$. Here the $e_j$ are distinct
and their sum is $0$. We denote
\begin{equation}
\omega_0=0,\quad\omega_1=1/2,\quad\omega_2=\tau/2,\quad\omega_3=(1+\tau)/2;
\end{equation}
then $e_k=\wp(\omega_k),\; 1\leq k\leq 3,$ and $\wp$ is periodic
with periods in $\Lambda_\tau$.

Let us define the functions $t(\tau)$ and $p(\tau)$ by
\begin{equation}\label{change}
t(\tau)=\frac{e_3(\tau)-e_1(\tau)}{e_2(\tau)-e_1(\tau)},\quad
y(t)=\frac{\wp(p(\tau)|\tau)-e_1(\tau)}{e_2(\tau)-e_1(\tau)}.
\end{equation}
The function $t(\tau)$ is the fundamental invariant of the group
$\Gamma[2]$ of the linear fractional transformations
represented by matrices $A\in SL(2,\Z)$ satisfying
$A\equiv I\; (\mod\; 2)$.

Then $p(\tau)$ satisfies the elliptic form of Painlev\'e VI,
\begin{equation}\label{elliptic}
\frac{d^2p(\tau)}{d\tau^2}=-\frac{1}{4\pi^2}\sum_{k=1}^3
\alpha_k\wp'(p(\tau)+\omega_k|\tau).
\end{equation}
Here
\begin{equation}\label{corres}
(\alpha_0,\alpha_1,\alpha_2,\alpha_3)=(\alpha,-\beta,\gamma,1/2-\delta).
\end{equation}
For the proof that (\ref{elliptic}) is equivalent to (\ref{1}) we refer to
\cite{Manin}.

Suppose that $y$ is an exceptional solution.
By (\ref{change}) this means
that
\begin{equation}\label{cond}
p(\tau)\not\equiv \omega_k\quad\mod\;\Lambda_\tau,\quad\tau\in H,
\quad k=0,\ldots,3.
\end{equation}
Moreover, as the only critical points of $z\mapsto\wp(z,\tau)$
are those congruent to $\omega_k$, we can locally solve the second equation
in (\ref{change}) with respect to $p$, and the implicit function theorem
implies that $p$ is holomorphic in $H$.

We use the following result of Earle \cite[Thm. 4.13]{Earle}:
\vspace{.1in}

\noindent
{\bf Theorem A.} {\em Let $p:H\to\C$ be
a holomorphic function with the property
that $p(\tau)\neq m+n\tau$ for all $\tau\in H$ and all integers $m,n$.
Then
\begin{equation}\label{p}
p(\tau)=\mu + \nu \tau ,
\end{equation}
where $\mu$ and $\nu$ are real, and $(\mu,\nu)\not\in\Z\times\Z$.}
\vspace{.1in}

Applying this theorem, we obtain that a solution $y(t)$ of (\ref{1})
described
by (\ref{change}) is exceptional
if and only if $p$ is of the form (\ref{p}),
with real $(\mu,\nu)\not\in(\Z/2)\times(\Z/2)$. Substituting to
(\ref{elliptic}), we obtain
\begin{equation}\label{reduction}
\sum_{k=0}^3\alpha_k\wp'(\mu + \nu \tau +\omega_k|\tau)\equiv 0.
\end{equation}
Such solutions are called {\em Picard's solutions}.
Picard \cite{Picard} found that they
exist in the case $\alpha_j=0$, $0\leq j\leq 3$. But of course they also
exist whenever (\ref{reduction}) holds. We mention the following
\vspace{.1in}

\noindent
{\bf Corollary of Theorem A.} {\em Let $y(t)$ be a multi-valued analytic function
in $\C\backslash\{0,1\}$, which has an analytic continuation along
every curve in $\C\backslash\{0,1\}$. Suppose that $y(t)\not\in\{0,1,t\}$
for all $t\in\C\backslash\{0,1\}$ and for all branches of $y$.
Then $y(t)$ is of the form (\ref{change}) with $p$ as in (\ref{p}). In
particular, $y$ is a solution of (\ref{1}) with
parameters $(0,0,0,1/2)$.}
\vspace{.1in}

Now our problem of classification of exceptional solutions is reduced to a
problem about elliptic functions: 

{\em For which
$\alpha_k,\mu,\nu$ do we have the identity (\ref{reduction})?}

To simplify (\ref{reduction}) we use the formulas
$$\wp(z+\omega_k)=e_k+\frac{(e_k-e_i)(e_k-e_j)}{\wp(z)-e_k}, \quad\{ i,j,k\}=\{1,2,3\}.$$
Differentiating these formulas with respect to $z$,
we obtain
$$\wp'(z+\omega_k)=-\frac{(e_k-e_i)(e_k-e_j)}{(\wp(z)-e_k)^2}\wp'(z),$$
and
substituting to (\ref{reduction}) we obtain after simplification
\begin{equation}\label{ma}
\alpha_0(w-e_1)^2(w-e_2)^2(w-e_3)^2=\sum_{k=1}^3
\alpha_k(w-e_i)^2
(w-e_j)^2(e_k-e_i)(e_k-e_j),
\end{equation}
where $w(\tau)=\wp(\mu + \nu \tau | \tau)$.
\vspace{.1in}

\noindent
{\bf Proposition 1.} {\em  If at least one $\alpha_k\neq 0$, the equation
(\ref{ma}) can only hold when $\mu,\nu$ are rational.}
\vspace{.1in}

{\em Proof.} The functions $e_j$ are the roots of the
equation
\begin{equation}\label{el}
4x^3-g_2x-g_3=0,
\end{equation}
whose coefficients are modular forms. In particular, if we set
$T\tau=\tau+1$, then the $g_j$ are invariant with respect to $T$
and thus the $e_j$ are invariant with respect to $T^3$. Then
it follows from (\ref{ma}) that $w(T^{18}\tau)=w(\tau)$.
Now
$$w(T^n\tau)=\wp(\mu+n\nu+\nu\tau|\tau+n)=\wp(\mu+n\nu+\nu\tau|\tau)=
\wp(\mu+\nu\tau+m+n\nu|\tau),$$
for all integers $m,n$. If $\nu$ is irrational we can arrange a sequence
$(m_k,n_k)$ such that $n_k$ are divisible by $18$, and $s_k=m_k+n_k\nu\to 0$.
Then
$$w(\tau)=w(\tau+s_k),\quad s_k\to 0,\quad s_k\neq 0,$$
which cannot happen for a non-constant analytic function.
This contradiction shows that $\nu$ is rational. A similar argument shows
that $\mu$ is also rational.
\vspace{.1in}

\noindent
{\bf Proposition 2.} { \em All exceptional solutions
of (\ref{1}) are algebraic.}
\vspace{.1in}

{\em Proof.} The function $w$ satisfying (\ref{ma}) can take only finitely many
values (at most 6) on any orbit of $\Gamma[2]$. Therefore $y$ can take only
finitely many values at each point. As $y$ omits $0,1,\infty$,
Picard's Great Theorem implies that the singularities at $0,1,\infty$ are
algebraic.
\vspace{.1in}

Actually one can write an explicit algebraic equation which all exceptional
solutions must satisfy. For this we express $w$ in terms of $y$ and the 
$e_j$ in terms of $t$ from
(\ref{change}) and substitute this expression to (\ref{ma}). We obtain:
\begin{equation}\label{polyeq}
\alpha_0y^2(y-1)^2(y-t)^2
-\alpha_1t(y-1)^2(y-t)^2-\alpha_2(1-t)y^2(y-t)^2
\end{equation}
$$-\alpha_3t(t-1)y^2(y-1)^2=0.$$

To determine how many exceptional solutions are possible,
one has to find for each $(\alpha_0,\ldots,\alpha_3)$
the number of irreducible factors of this equation,
and to check which of these factors define
algebraic solutions of (\ref{1}).
For example, in the case considered in \cite{Lin},
when all $\alpha_j$ are equal, we obtain three factors:
$$ y^2(y-1)^2(y-t)^2-t(y-1)^2(y-t)^2-(1-t)y^2(y-t)^2+t(t-1)y^2(y-1)^2$$
\begin{equation}\label{lin2}
=(y^2-t)(y^2-2y+t)(y^2-2yt+t).
\end{equation}
In this case, each of the three factors determines a solution of (\ref{1}). 
So we obtained a simple proof of the main theorem of \cite{Lin}.
We can state the result as follows:
\vspace{.1in}

\noindent
{\bf Proposition 3.} {\em When not all $\alpha_j=0$, exceptional solutions
are algebraic functions given by the polynomial equation (\ref{polyeq}).
Their number is the number of non-trivial irreducible
factors of this polynomial that satisfy (\ref{1}).
A factor is called non-trivial if it depends on both $y$ and $t$
and is not a constant multiple of $y-t$.}
\vspace{.1in}

It is easy to see that an exceptional solution cannot be rational,
see, for example, \cite[Proposition 6]{Douady},
so the number of exceptional solutions
is at most $3$, and they are at most $6$-valued. 

Next we determine 
all cases when the polynomial in
(\ref{polyeq}) is reducible.
\vspace{.1in}

\noindent
{\bf Proposition 4.} {\em When $\alpha_3=0$, the polynomial
(\ref{polyeq}) factors as
$$(y-t)^2P_0(y,t),$$
where
$$P_0(y,t)=\alpha_0(y-1)^2y^2-\alpha_2y^2-
t\left(\alpha_1(y-1)^2-\alpha_2y^2\right)$$
has no non-trivial factors.}
\vspace{.1in}

In this case we may have at most one  exceptional
solution defined by $P_0(y,t)=0$.
\vspace{.1in}

\noindent
{\bf Proposition 5.} {\em If $\alpha_3\neq 0$, then the polynomial
(\ref{polyeq}) has a non-trivial
factorization if 
$$\alpha_j=u_j^2,\quad\mbox{where}\quad\sum_{j=0}^3\pm u_j=0,$$
for any choice of signs. An equivalent condition is
$$\left( \alpha_0^2+\alpha_1^2+\alpha_2^2 +\alpha_3^2-2(\alpha_0\alpha_1+
\alpha_0\alpha_2+\alpha_0\alpha_3+\alpha_1\alpha_2+\alpha_1\alpha_3+
\alpha_2\alpha_3)\right)^2$$
\begin{equation}\label{condition}
=64\alpha_0\alpha_1\alpha_2\alpha_3.
\end{equation}
The surface defined by (\ref{condition}) contains three lines
\begin{equation}
\label{line1}
\left\{\alpha_0=\alpha_1,\;\alpha_2=\alpha_3\right\},\end{equation}
\begin{equation}
\left\{\alpha_0=\alpha_2,\;\alpha_1=\alpha_3\right\},\label{line2}
\end{equation}
\begin{equation}\label{line3}
\left\{\alpha_0=\alpha_3,\;\alpha_1=\alpha_2\right\}.
\end{equation}
The polynomial (\ref{polyeq}) is a product of three non-trivial ir\-re\-du\-ci\-ble
fac\-tors if
\newline
 $(\alpha_0,\alpha_1,\alpha_2,\alpha_3)$ belongs
to one of these lines, and $\alpha_0\alpha_1\alpha_2\alpha_3\neq 0$.
}
\vspace{.1in}

Using Maple and Mathematica we determined that the cases listed in 
Proposition 5 exhaust all factorizations of our polynomial.
But this fact will not be used in the proof of our main result.

When (\ref{condition}) does not hold, computation indicates that
the polynomial (\ref{polyeq}) is
irreducible, and can define at most one exceptional solution of (\ref{1}).
It remains to determine which algebraic functions defined by factors of
(\ref{polyeq})
actually satisfy (\ref{1}).
Computation with Maple indicates that in the case when (\ref{polyeq})
is irreducible, the resulting algebraic function with 6 branches
does not satisfy (\ref{1}). 
In the case (\ref{condition}) when (\ref{polyeq}) splits into two
irreducible factors, algebraic functions arising from these factors are of
degrees $2$ and $4$. The function determined by the factor of
degree $2$ never satisfies (\ref{1}), while the function determined
by the factor of degree $4$ satisfies (\ref{1}) if and only if 
three of the $\alpha_j$ are equal and the fourth is equal to the
sum of these three. 
If one of the equations (\ref{line1}), (\ref{line2}), or (\ref{line3})
is satisfied, then one of the factors satisfies the equation
and the other two factors do not.

We will prove all these facts below without reliance on a computer.
Our main result is the following.
\vspace{.1in}

\noindent
{\bf Theorem 1.} {\em The complete list of
exceptional solutions of Painlev\'e VI is 
the following:

\noindent
If $\alpha_j=0,\; 0\leq j\leq 3$, they are Picard's solutions
with real $(\mu,\nu)\in \R^2\backslash\Z^2$.

\noindent
If $\alpha_0=\alpha_1,\quad\alpha_2=\alpha_3,$ then there is a solution
\begin{equation}\label{a}
y(t)=\sqrt{t}.
\end{equation}
If $\alpha_0=\alpha_2,\quad\alpha_1=\alpha_3,$ then there is a solution
\begin{equation}\label{b}
y(t)=1+\sqrt{1-t}.
\end{equation}
If $\alpha_0=\alpha_3,\quad\alpha_1=\alpha_2,$ then there is a solution
\begin{equation}\label{c}
y(t)=t+\sqrt{t^2-t}.
\end{equation}
If $\alpha_0=9\alpha_1=9\alpha_2=9\alpha_3\neq 0,$ then there is
a unique  solution defined by
\begin{equation}\label{d}
3y^4-4ty^3-4y^3+6ty^2-t^2=0.
\end{equation}
If $9\alpha_0=\alpha_1=9\alpha_2=9\alpha_3\neq 0,$ then there is a
unique solution
defined by
\begin{equation}\label{e}
y^4-6ty^2+4t(t+1)y-3t^2.
\end{equation}
If $9\alpha_0=9\alpha_1=\alpha_2=9\alpha_3\neq 0,$ then there is
a unique solution defined by
\begin{equation}\label{f}
y^4-4y^3+6ty^2-4t^2y+t^2. 
\end{equation}
If $9\alpha_0=9\alpha_1=9\alpha_2=\alpha_3\neq 0,$ then there
is a unique solution defined by
\begin{equation}\label{g}
y^4-4ty^3+6ty^2-4ty+t^2.
\end{equation}}
\vspace{.1in}

Equations (\ref{a}), (\ref{b}), (\ref{c}) are permuted
by the group generated by
\begin{equation}\label{group}
(t,y)\mapsto(1-t,1-y),\quad\mbox{and}\quad(t,y)\mapsto(1/t,y/t),
\end{equation}
which is isomorphic to $S_3$. Equations (\ref{d}), (\ref{e}), (\ref{f}),
(\ref{g}) are permuted by the group isomorphic to $S_4$ which is
obtained by adding the transformation
\begin{equation}\label{group2}
(t,y)\mapsto(1/t,1/y)
\end{equation}
to (\ref{group}).

All curves (\ref{d}), (\ref{e}), (\ref{f}), (\ref{g}) are of genus zero.
A uniformization of (\ref{d}) is
$$y=\frac{1}{1-z^2},\quad t=\frac{2z-1}{(z-1)^3(z+1)}.$$
The rest are obtained by substitutions (\ref{group2}), (\ref{group}):
$$y=1-z^2,\quad t=\frac{(z+1)(z-1)^3}{2z-1} \quad\quad\mbox{for}\quad (\ref{e}),$$
$$y=z^2,\quad t=-\frac{z^3(z-2)}{2z-1} \quad\quad\mbox{for}\quad (\ref{f}),$$
$$y=-\frac{2z-1}{z(z-2)},\quad t=-\frac{2z-1}{z^3(z-2)}\quad\quad\mbox{for}\quad
(\ref{g}).$$

{\em Proof of Theorem 1.} We have proved that all exceptional solutions
are Picard solutions parametrized by
$$y(\tau)=\frac{\wp(\mu +\nu\tau|\tau)-e_1(\tau)}{
e_2(\tau)-e_1(\tau)}, \quad t=\frac{e_3(\tau)-e_1(\tau)}{e_2(\tau)-e_1(\tau)}.$$
Two such solutions are the same (obtained by an analytic continuation)
if, and
only if, 
$$(\mu,\nu)A=\pm(\mu,\nu)
\quad (\mod\; \Z\times\Z),\quad\mbox{where}\quad A\in\Gamma[2].$$
Let us say that two rational vectors $(\mu,\nu)$ and $(\mu',\nu')$
are equivalent if
$$(\mu,\nu)=\pm(\mu',\nu')\quad (\mod\;\Z\times\Z).$$
Then  the group $\Gamma[2]$ acts on the equivalence classes,
and we need the list of all classes whose orbit has length at most $6$.

We have the following
\vspace{.1in}

\noindent
{\bf Lemma 1.} \cite[Lemma 3]{Mazz} {\em Every $\Gamma[2]$ orbit
contains a vector equivalent to one of the
following: 
\begin{equation}\label{cases}
(0,M/N),\quad (M/N,0),\quad (M/N,M/N),
\end{equation}
where $M,N$ are defined as follows. Let $\mu=\mu_1/\mu_0$, $\nu=\nu_1/\nu_0$
be the reduced representations. Then $N$ is the least common multiple
of $\mu_0,\nu_0$ and $M$ is the greatest common divisor of
$\mu_1N/\mu_0,\;\nu_1N/\nu_0$, so that $M,N$ are coprime, and
$$\mu=mM/N,\quad\nu=nM/N,$$
where $m,n$ are coprime.

Then the orbit of $(\mu,\nu)$ contains a vector
of the list (\ref{cases})
if, and only if, $(m,n)$ is
(even,odd), (odd,even) or (odd,odd), respectively.

The three solutions corresponding to the vectors (\ref{cases}) are permuted
by the group generated by (\ref{group})
}
\vspace{.1in}

To understand the orbits completely, it remains to check which of the
three vectors (\ref{cases}) are on the same orbit.
\vspace{.1in}

\noindent
{\bf Lemma 2.} {\em Let $M,N$ be coprime integers. If $N$ is odd, then
the three points (\ref{cases}) are in one $\Gamma[2]$ orbit.
If $N$ is even, they are in three distinct orbits}.
\vspace{.1in}

{\em Proof.} The vectors $(0,M/N)$ and $(M/N,0)$ are on the same orbit if, and only if, there
exists a matrix in $\Gamma[2]$ such that
\begin{equation}\label{eq}
\left(\begin{array}{cc}a&b\\ c&d\end{array}\right)\left(\begin{array}{c}M/N\\
0\end{array}\right)=\pm\left(\begin{array}{c}0\\ M/N\end{array}\right)\quad
(\mod\; \Z\times\Z),
\end{equation}
which is equivalent to
\begin{equation}
\label{11}
a\equiv 0\quad (\mod\; N),
\end{equation}
and
\begin{equation}\label{2}
cM\equiv\pm M\quad (\mod\; N).
\end{equation}
As $a$ is odd, we conclude from (\ref{11}) that $N$ must be odd.

In the opposite direction, if $N$ is odd, we can take
\begin{equation}\label{mat}
\left(\begin{array}{cc}a&b\\ c&d\end{array}\right)=
\left(\begin{array}{cc}-N&N+1\\ -1-N^2& 1+N(N+1)\end{array}\right)\in\Gamma[2],
\end{equation}
and (\ref{eq}) will be satisfied. 

Now let us investigate when $(M/N,0)$ and $(M/N,M/N)$
are on the same orbit. We have
\begin{equation}\label{eq2}
\left(\begin{array}{cc}a&b\\ c&d\end{array}\right)\left(\begin{array}{c}M/N\\
M/N\end{array}\right)=\pm\left(\begin{array}{c}M/N\\ 0\end{array}\right)\quad
(\mod\; \Z\times\Z),
\end{equation}
which is equivalent to
\begin{equation}\label{3}
c+d\equiv 0\quad(\mod\; N),
\end{equation}
and
\begin{equation}\label{4}
(a+b)M=\pm M\quad(\mod\; N).
\end{equation}
As $c+d$ is always odd, we conclude from (\ref{3}) that $N$ must be odd.
In the opposite direction,
if $N$ is odd, use the same matrix as in (\ref{mat})
and (\ref{eq2}) will be satisfied. 

Now it is easy to find the number of elements in an orbit.
When $N$ is odd, all vectors $(\mu_1/N,\nu_1/N)$
with the greatest common factor of $\mu_1,\nu_1$ coprime to $N$
belong to one orbit. This orbit is of
length greater than $6$ when $N\geq 5$.

When $N$ is even, such vectors lie on three orbits of equal length,
corresponding to the three vectors (\ref{cases}), when $N=4$ we have three
orbits of length $2$, and when $N=6$ we have three orbits of length $4$.

Thus the result is that exceptional solutions correspond to the following pairs
$(\mu,\nu)$:

a) $(1/4,0),\; (0,1/4),\; (1/4,1/4)$, representing
three distinct two-valued solutions.

b) $(1/3,1/3)$, representing one four-valued solution.

c) $(1/6,0),\; (0,1/6),\; (1/6,1/6)$, representing three four-valued solutions.

As one equation cannot have two different four-valued exceptional solutions,
the three solutions in c) must belong to different equations. The three
solutions in a) may belong to one equation, or to three different equations.
The group (\ref{group}) permutes solutions of the type a)
and permutes solutions
of the type b).

One can verify that vectors a) correspond to solutions (\ref{a}), (\ref{b}),
(\ref{c}) in Theorem 1, vector b) corresponds to (\ref{d}),
and vectors, c) correspond to the remaining three
solutions (\ref{e}), (\ref{f}), (\ref{g}).

For a) this is easy. To check b) and c), we write the tripling formula
for the elliptic function
$$w(z)=\frac{\wp(z)-e_1}{e_2-e_1},$$
which can be obtained from the well-known addition theorem for $\wp$.
We have
\begin{equation}\label{tripling}
w(3z)=y\left(\frac{y^4+4yt-6y^2t-3t^2+4yt^2}{4y^3t-6y^2t+4y^3-3y^4+t^2}
\right)^2=:y\left(\frac{f(y,t)}{g(y,t)}\right)^2,
\end{equation}
where $y=w(z)$.
At the points $z$ of third order, we have $w(3z)=\infty$, while at the
points of $6$-th order, $w(3z)\in\{0,1,t\}$. So we have to solve four equations
$$f(y,t)=0,\quad g(y,t)=0,\quad yf(y,t)-g(y,t)=0,\quad yf(y,t)-tg(y,t)=0.$$
The first two polynomials are
irreducible. Factoring the other two
we obtain:
$$f(y,t)-g(y,t)=(y-1)(4y^3-y^4-6y^2t-t^2+4yt^2)^2$$
and
$$f(y,t)-tg(y,t)=(y-t)(y^4-4yt+6y^2t-4y^3t+t^2)^2,$$
which together with $f$ and $g$ gives the four polynomials
in (\ref{d}), (\ref{e}), (\ref{f}), (\ref{g}). 
\vspace{.1in}

{\bf Remarks.}
\vspace{.1in}

1. The equation
(\ref{condition}) defines a Kummer surface \cite[p. 21, footnote]{Hu}.
\vspace{.1in}

2. All algebraic solutions of Painlev\'e VI have been classified in
\cite{LT}. However this classification is only up to
B\"acklund transformations, and B\"acklund transformations in general do not
map exceptional solutions to exceptional solutions \cite{Lin}.
\vspace{.1in}

We thank Chang-Shou Lin who sent us \cite{Lin} and
communicated the solution (\ref{lin})
before \cite{Lin2} was released, and Clifford Earle for showing us
an unpublished letter of Carleson with an elementary proof of Theorem A.
We also thank Oleg Lisovyy, T. N. Venkataramana
and Math Overflow participants \cite{MO} for illuminating
discussions.

{\em A. E. and A. G.: Department of Mathematics, Purdue University,
West Lafayette Indiana 47907,

A. H.: Department of Mathematics, University of Illinois, Urbana,
Illinois 61801.}

\end{document}